\input amstex
\documentstyle{amsppt}
\input bull-ppt
\keyedby{bull269e/jxs}


\define\csus{\notSigma}  
\define\bbq{\Bbb Q}   
\define\equdef{\buildrel \operatorname{def} \over =}


\topmatter
\cvol{26}
\cvolyear{1992}
\cmonth{April}
\cyear{1992}
\cvolno{2}
\cpgs{264-268}

\title A theory of algebraic cocycles\endtitle
\author Eric M. Friedlander and H. Blaine Lawson, 
Jr.\endauthor
\shortauthor{E. M. Friedlander and H. B. Lawson, Jr.}
\date April 4, 1991\enddate
\thanks The first author was partially supported by the 
IHES, NSF~Grant
DMS 8800657, and NSA~Grant MDA904-90-H-4006.  The second 
author was 
partially supported by NSF~DMS 8602645\endthanks
\subjclass Primary 14F99, 14C05\endsubjclass
\keywords Algebraic cycle, Chow variety, algebraic 
cocycle, co\-homology\endkeywords
\address Department of Mathematics, Northwestern 
University, Evanston, Illinois
60208\endaddress
\address Department of Mathematics, State University of 
New York,
Stony Brook, New York 11794\endaddress
\abstract 
We introduce the notion of an algebraic cocycle as the 
algebraic
analogue of a map to an Eilenberg-MacLane space.  Using 
these cocycles we
develop a  ``cohomology theory" for complex algebraic 
varieties.  The theory is
 bigraded, functorial, and admits Gysin maps.  It carries 
a natural cup
product and a pairing to $L$-homology.  Chern classes of 
algebraic bundles are
defined in the theory.  There is a natural transformation to
(singular) integral cohomology theory that preserves cup 
products.
Computations in special cases are carried out.  On a
smooth variety it is proved that there are algebraic 
cocycles in each
algebraic rational $(p,p)$-cohomology class.  \endabstract

\endtopmatter


\document

In this announcement we present the outlines of a cohomology
theory for  algebraic varieties based on a
new concept of an algebraic cocycle.
Details will appear in [FL]. Our cohomology is a companion
 to the $L$-homology theory recently studied in [L,
F, L-F1, L-F2, FM].  This homology is a bigraded theory
based directly on the structure of the space of algebraic 
cycles.  
It admits a natural  transformation to  integral
homology that generalizes the usual map taking a cycle to 
its
homology class.  Our new cohomology theory is similarly
bigraded and based on the structure of the space of 
algebraic
cocycles.  It carries a ring structure coming from the 
complex
join (an elementary construction of projective geometry), 
and it
admits a natural transformation $\Phi$ to  integral
cohomology.  Chern classes are  defined in  the theory
and transform under $\Phi$  to the usual ones.  Our 
definition of
cohomology is very far from a duality construction on
$L$-homology.  Nonetheless, there is a natural and 
geometrically
defined Kronecker pairing between our ``morphic 
cohomology'' and $L$-homology.
\par \rm  The foundation stone of our theory is the notion 
of an
effective algebraic cocycle, which is of some independent
interest.  Roughly speaking, such a cocycle on a variety 
$X$, with
values in a projective variety $Y$, is a morphism from $X$ 
to the
space of cycles on $Y$.  When $X$ is normal, this is 
equivalent
(by ``graphing") to a cycle on $X\times Y$ with 
equidimensional
fibres over $X$.  Such cocycles abound in algebraic 
geometry and
arise  naturally in many circumstances.  The simplest 
perhaps is
that of a flat morphism $f\:X \rightarrow Y$ whose
corresponding  cocycle associates to $x \in X$, the pullback
cycle $f^{-1}(\{ x \})$.  Many more arise naturally from 
synthetic
constructions in projective geometry.  We  show that every
variety is rich in cocycles.  Indeed if $X$ is  smooth  and
projective, then every rational cohomology class that is
Poincar\' e dual to an algebraic cycle is represented by
(i.e., is $\Phi \otimes\bbq$ of)  an algebraic cocycle.
\par\rm
In what follows the word {\it variety} will denote a 
reduced,
irreducible,  locally closed subscheme of some complex 
projective
space.  Such a variety is called projective if
it is in fact closed (equivalently, compact)
in some projective space.

\dfn{Definition 1} 
Given a projective variety $Y \subset
  P  \sp{ n}$  with a fixed embedding, we denote by
$ C \sb{ d} \sp{ s}(Y)$ the algebraic set (the  ``Chow 
set")  of
effective cycles of codimension-\<$s$ and degree $d$ with 
support
in $Y$.  For an arbitrary variety $X$ we then define an 
{\it effective algebraic cocycle\/} 
on $X$ with  values in $Y$ to be a
continuous algebraic map $\varphi \: X \rightarrow  C \sb{ 
d}
\sp{ s}(Y)$ \ (i.e., a morphism 
$\varphi \: \tilde X \rightarrow  C \sb{ d} \sp{ s}(Y)$
from the weak normalization $\tilde X$ of $X$).  The space 
of all
such cocycles provided with the compact-open topology will 
be
denoted by
$ C \sb{ d} \sp{ s}(X;Y)$.
\enddfn


Note that $ C \sb{ d} \sp{ s}(X;Y)$ is {\it a priori}  a 
hybrid
construction consisting of  algebrogeometric objects but
carrying the compact-open topology. However, for normal
projective varieties $X$, it  has a purely algebraic 
description. 

\thm{Proposition 2} 
If
$     X$ is  projective and normal, then the space
$ C \sb{ d} \sp{ s}(X;Y)$ admits the structure
of a locally closed, reduced subscheme of some complex
projective space. 
\ethm

Formal direct sum determines  
 an abelian topological monoid structure on the spaces
$$
C \sp{ s}(X;Y) \equdef  \coprod  \sb{ d \geq
0}C \sb{ d} \sp{ s}(X;Y)\.
$$
This structure is proved to be independent of the projective
embedding chosen for $Y$. In analogy with the construction 
of 
 $L$-homology, we have the following definition. 

\dfn{Definition 3} 
Let
$     X$ and
$     Y$ be varieties, with
$     Y$ projective.  Denote by
$ Z \sp{ s}\!(X\!;\!Y\!)$ the homotopy theoretic
group completion of
$ C \sp{ s}(X;Y)$ (i.e., the loops on the
classifying space of
$ C \sp{ s}(X;Y))$,
$$
Z \sp{ s}(X;Y) \equiv  \Omega   B(C \sp{
s}(X;Y))\.
$$
Then the  
{\it bivariant morphic cohomology\/} 
of
$     X$ with coefficients in
$     Y$ is defined to be the homotopy groups
of
$ Z \sp{ s}(X;Y), $
$$
L \sp{ s}H \sp{ q}(X;Y) \equiv  \pi  \sb{
2s-q}(Z \sp{ s}(X;Y)),\qquad 2s \geq q \geq
0\.
$$
\enddfn

The first fundamental result concerning these spaces is the
Algebraic Suspension Theorem,  which asserts that the
algebraic suspension maps 
$C \sb{ d} \sp{ s}(X;Y) \rightarrow
C \sb{ d} \sp{ s}(X;\csus Y)$,  introduced in [L], induce 
a homotopy
equivalence
$$
Z \sp{ s}(X;Y) \ \buildrel \cong \over \longrightarrow 
\  Z \sp{ s}(X;\csus Y), 
$$
and thus an isomorphism
$$
L \sp{ s}H \sp{ q}(X;Y) \ \cong \ L \sp{ s}H \sp{ q}(X; 
\csus Y)
$$
for all $s$ and $q$.
\par\rm
Although this theory has been
developed in the ``bivariant context" of Definition 3, we 
shall
focus our attention here on the important special case in 
which
$     Y$ is some projective space
$  P  \sp{ N}. $
We know from [L] that
$  \Omega   B(C \sp{ s}( P  \sp{ N}))$
has the homotopy type of a generalized Eilenberg-MacLane
space
$$
 \Omega   B(C \sp{ s}( P  \sp{
N})) \cong K( {\bold Z} ,0) \times K( {\bold Z}
,2) \times \cdots \times K( {\bold Z} , 2s),
\qquad N \geq s\.
$$
In particular, this homotopy type is independent
of
$ N \geq s.$  We conclude that
$  \Omega   B(C \sp{ s}( P  \sp{ N}))$
``modulo"
$  \Omega   B(C \sp{ s-1}( P  \sp{
N-1}))$ represents cohomology.  This motivates
the following

\dfn{Definition 4} 
For any variety
$     X$ and any
$ s \geq 0,$ let
$ Z \sp{ s}(X)$ denote the homotopy fibre
of the natural map
$ B  C \sp{ s-1}(X, P  \sp{ s-1}) \rightarrow
B  C \sp{ s}(X, P  \sp{ s}) $,
$$
Z \sp{ s}(X) \equiv \operatorname{htyfib} \{ B 
C \sp{ s-1}(X, P  \sp{ s-1}) \rightarrow
B  C \sp{ s}(X, P  \sp{ s}) \}\.
$$
For any
$ 0 \leq q \leq 2s,$ the  
{\it morphic cohomology\/}
group
$ L \sp{ s}H \sp{ q}(X)$ is defined by
$$
L \sp{ s}H \sp{ q}(X) \equiv  \pi  \sb{
2s-q}(Z \sp{ s}(X))\.
$$
\enddfn

Definitions 3 and 4 are related as follows.

\thm{Theorem 5} 
For any variety
$     X$ and any
$ N \geq s \geq 0,$ there is a natural homotopy
equivalence 
$$
Z \sp{ s}(X; P  \sp{ N}) \cong Z \sp{
0}(X) \times Z \sp{ 1}(X) \times \cdots \times
Z \sp{ s}(X)\.
$$
\ethm

The splitting asserted in Theorem 5 arises
from natural maps
$$
SP \sp{  \infty }( P  \sp{ N}) \rightarrow
SP \sp{  \infty }( P  \sp{ k}),\qquad k \leq
N
$$
obtained by viewing
$  P  \sp{ N}$ as
$ SP \sp{ N}( P  \sp{ 1})$ and
$  P  \sp{ k}$ as
$ SP \sp{ k}( P  \sp{ 1}),$ where
$ SP \sp{ j}( P  \sp{ n})$ denotes the j-fold
symmetric product of
$  P  \sp{ n}. $
\par\rm
We view morphic cohomology as the theory
corresponding to ``algebraic" as opposed
to ``arbitrary continuous" maps from
$     X$ into Eilenberg-MacLane spaces.
This perspective is formalized in the following.

\thm{Theorem 6} 
For any variety $X$, the elementary
complex join operation induces a natural ring structure on
 $$
L \sp{  \cdot }H \sp{  \cdot }(X) \equiv
 \bigoplus  \sb{ s \geq 0}L \sp{ s}H \sp{
\cdot }(X)\.
$$
 Furthermore,
there is a natural transformation of graded
rings
$$
 \Phi \: L \sp{  \cdot }H \sp{  \cdot }(X)
 \rightarrow H \sp{ \cdot }(X; {\bold Z} ),
$$
and, if
$     X$ is projective then
$$
\Phi (L \sp{ s}H \sp{ 2s-j}(X)) \otimes
{\bold C}\,  \subset \,H \sp{ s,s-j} \oplus
H \sp{ s-1,s-j+1} \oplus \cdots \oplus H \sp{
s-j,j },
$$
where
$ H \sp{ p,q}$ denotes the
$ (p,q)$\<th Dolbeault component of
$ H \sp{ p+q}(X; {\bold C} ). $
\ethm

The existence of a ring structure in morphic
cohomology provides it with a structure not
possessed by $L$-homology.  On the other
hand, the natural operations on $L$-homology
constructed in [FM] via the join
operation naturally determine operations
on our morphic cohomology groups.
\par\rm
The restriction on the image of
$  \Phi $ given in Theorem 6 is complemented
by the following existence result.  A key ingredient in 
its proof 
is the total Chern class map of [LM],
$$
BU \sb{ s} \rightarrow C \sp{ s}( P  \sp{ \infty }),
$$
which can be viewed geometrically as the inclusion of 
degree 1
cycles into the space of all cycles on 
$ P \sp{ N}$ for N sufficiently large.

\thm{Theorem 7} 
Let
$     E$ be a vector bundle over $X$
 generated by its global
sections, over a variety
$     X$.  Then there are naturally defined
chern classes
$$
c \sb{ k}(E) \in L \sp{ k}H \sp{ 2k}(X)
$$
with the property that
$  \Phi (c \sb{ k}(E)) \in H \sp{ 2k}(X; {\bold Z} )$
is the usual $k$\<th chern class of
$     E$.  Consequently, if
$     X$ is a smooth projective variety,
then the Poincare dual of the fundamental
class of each algebraic subvariety lies in the
subring of
$ H \sp{ *}(X; {\bold Z} )$ generated
$  \Phi (L \sp{  \cdot }H \sp{  \cdot }(X)) $. 
\ethm

Not surprisingly, codimension-1 morphic cohomology
is the easiest to compute.  We have the following
computation.

\thm{Theorem 8}  
Let $ X$ be a projective variety.  Then
$$
L^1H^q(X)=\cases \bold Z&\text{if }q=0,\\
H^1(X;\bold Z)&\text{if }q=1,\\
H^{1,1}(X;\bold Z)&\text{if }q=2,\\
0&\text{otherwise},\endcases
$$
where
$ H \sp{ 1,1}(X; {\bold Z} )$ denotes
$ H \sp{ 2}(X; {\bold Z} ) \cap  \rho  \sp{
-1}H \sp{ 1,1}(X; {\bold C} )$  and where
$  \rho \:
 H \sp{ 2}(X; {\bold Z} ) \rightarrow H \sp{
2}(X; {\bold C} )$ is the coefficient homomorphism. 
\ethm

In [F] a result similar to Theorem 8 but
only applying to smooth projective varieties
was proved for $L$-homology.  This suggests
that morphic cohomology and $L$-homology
should satsify some form of duality.  One
possible candidate for a possible duality
pairing is given in the following proposition.

\thm{Proposition 9} 
For any variety
$     X,$ there is a natural Kronecker pairing
between $L$-homology and morphic cohomology,
$$
L \sp{ s}H \sp{ q}(X) \otimes L \sb{ r}H \sb{
q}(X) \rightarrow  {\bold Z},
$$
which is naturally compatible with the usual
Kronecker pairing  
$$
  H \sp{ q}(X; {\bold Z} ) \otimes H \sb{
q}(X; {\bold Z} ) \rightarrow  {\bold Z} \.
$$
\ethm

The definition of this pairing is pleasingly
geometric.  Namely, given a cycle
$     W$ on
$ X \times  P  \sp{ N}$ and a cycle
$     Z$ on X, we take the image in
$ H \sb{ *}( P  \sp{ N}; {\bold Z} )$ of the
fundamental class of the restriction of
$     W$ to
$     Z. $


\enddocument